\documentclass{amsart} 
\usepackage{amssymb, amscd}
\numberwithin{equation}{section}

\def\CC{{\mathbb C}}  

\def\EE{{\mathbb E}}    
\def\FF{{\mathbb F}}  
\def\GG{{\mathbb G}}

\def\PP{{\mathbb P}}
\def\QQ{{\mathbb Q}}

\def\ZZ{{\mathbb Z}}

\def\Acal{{\mathcal A}}

\def\Fcal{{\mathcal F}}

\def\Xcal{{\mathcal X}}

\let\co\colon
\newcommand\Definition[1]{\emph{#1}}
\newcommand\map[3]{\ensuremath{{#1}\co{#2}\to{#3}}}

\newcommand\proofsquare{\nobreak\hfill \hbox{%
\vrule height 5pt 
\kern-.4pt
 \vbox{%
\hrule width 5pt depth0pt height.4pt
 \kern4.6pt \hrule  }
\kern-3.75pt 
\vrule height 5pt}\kern1pt
\par}

\newtheorem{theorem}{Theorem}[section]

\newtheorem{proposition}[theorem]{Proposition}

\newtheorem{definition-lemma}[theorem]{Definition-Lemma}

\theoremstyle{definition}

\newtheorem{example}[theorem]{Example}

\theoremstyle{remark} 
\newtheorem{remark}[theorem]{Remark}

\begin{document}

\title[Cycles Representing The Top Chern Class]{Cycles Representing The Top
Chern Class \\
of the Hodge Bundle on the Moduli Space\\
of Abelian Varieties
} 
\author{Torsten Ekedahl}
\address{Department of Mathematics\\
 Stockholm University\\
 SE-106 91  Stockholm\\
Sweden}
\email{teke@math.su.se}
\author{Gerard van der Geer}

\address{Faculteit Wiskunde en Informatica, University of
Amsterdam, Plantage Muidergracht 24, 1018 TV Amsterdam, The Netherlands}

\email{geer@science.uva.nl} 


\begin{abstract}

We give a generalization to higher genera of the famous formula $12 \,
\lambda=\delta$ for genus $1$. We also compute the classes of certain 
strata in the Satake compactification as elements of the push down of 
the tautological ring.

\end{abstract}

\maketitle

\begin{section}{Introduction}
\label{sec: intro}
\bigskip
\noindent
The fact that there exists a cusp form of weight $12$ on $\mathrm{SL}(2,\ZZ)$
with a simple zero in the cusp and no zero on the upper half plane
translates into the cycle relation $12 \lambda=\delta$. 
Here $\lambda$ is the divisor class corresponding
to to the factor of automorphy 
$((\begin{matrix}a & b\\ c& d\\  \end{matrix}),z) \mapsto (cz+d)$
and $\delta$ represents the class of the cusp. We wish to generalize this
relation to a relation in the Chow ring of the moduli of principally
polarized abelian varieties. The analogue of the class $\lambda$ is the top
Chern class $\lambda_g$ of the Hodge bundle and the analogue of $\delta$ 
is a codimension $g$ class $\delta_g$ living on
the boundary of the moduli space. The analogue
of the formula is then a relation of the form
$$
\lambda_g= (-1)^g \zeta(1-2g)\, \delta_g,
$$
with $\zeta(s)$ the Riemann zeta function. We now formulate a precise
version of this.

Let ${\Acal}_g/ \ZZ$ denote the moduli stack of principally polarized abelian
varieties of dimension~$g$. This is an irreducible algebraic stack of relative
dimension $g(g+1)/2$ with irreducible fibres over $\ZZ$. The stack ${\Acal}_g$
carries a locally free sheaf $\EE$ of rank $g$, the Hodge bundle, defined as
follows.  If $A/S$ is an abelian scheme over $S$ with zero section $s$ we get a
locally free sheaf $s^*\Omega^1_{A/S}$ of rank $g$ on $S$ and this is compatible
with pull backs. If $\pi: A \to S$ denotes the structure map it satisfies the
property $\Omega^1_{ A/S}= \pi^*(\EE)$ and its top 
Chern class $\lambda_g(A/S):=c_g(\Omega^1_{ A/S})$ 
(in the Chow ring of $S$) is then 
a pullback of a corresponding class in the universal case $\lambda_g :=
c_g(\EE)$. The Hodge bundle can be extended to a locally free sheaf (again
denoted by) $\EE$ on every smooth toroidal compactification $\tilde {\Acal}_g$
of ${\Acal}_g$ of the type constructed in \cite{F-C}, see Ch.\ VI,4 there. By a
slight abuse of notation we will continue to use the notation $\lambda_g$ for
its top Chern class.

The class  $\lambda_g$ is defined over $\ZZ$ and gives for each fibre
${\Acal}_g\otimes k$ with $k$ a field 
rise to a class, also denoted $\lambda_g$, 
in the Chow group $CH^g({\Acal}_g\otimes k )$, and in 
$CH^g(\tilde{\Acal}_g\otimes k)$. 
It was proved  in \cite{vdG1} that $\lambda_g$ vanishes in the Chow group
$CH^g_{\QQ}({\Acal}_g)$ with rational coefficients; however, it does not
vanish on ${\Acal}_g\otimes k$  and  in \cite{E-vdG} we studied the order 
of the torsion class $\lambda_g$. It also does not vanish in the Chow
group $CH^g_{\QQ}(\tilde{\Acal}_g\otimes k)$. Therefore one may ask for an
effective cycle representing the class $\lambda_g$ on a compactification.

There are several compactifications of ${\Acal}_g$. 
We let ${\Acal}_g^*$ be the minimal or Satake
compactification as defined in \cite{F-C}. This compactification
${\Acal}_g^*$ is a disjoint union
$$
{\Acal}_g^* = {\Acal}_g \sqcup {\Acal}_{g-1}\sqcup \ldots \sqcup {\Acal}_0.
$$

If $\tilde{\Acal}_g$ is a
suitable smooth toroidal compactification $\tilde{\Acal}_g$ as constructed in
\cite{F-C} we have a natural map $q\colon \tilde{\Acal}_g \to {\Acal}_g^*$ to
the Satake compactification. 
The moduli space ${\Acal}_g^{\prime}$ 
of rank $1$ degenerations is by definition
the inverse image of ${\Acal}_g \sqcup
{\Acal}_{g-1} \subset {\Acal}_g^*$ under the natural map 
$ q\colon \tilde{\Acal}_g \to {\Acal}_g^*$. 
The important fact is that
the space ${\Acal}_g'$ does {\sl not} depend on a choice $\tilde {\Acal}_g$ of
compactification of ${\Acal}_g$; it  is a \emph{canonical}
partial compactification
on ${\Acal}_g$. If we want a full  compactification then there is not really
a unique one, but we must make choices, see \cite{M2}.

The space ${\Acal}_g^{\prime}$ parametrizes semi-abelian varieties
with torus rank $\leq 1$. We let $\Delta_g$ be the closed locus of
${\Acal}_g^{\prime}$ which parametrizes the semi-stable abelian
varieties which are trivial extensions
$$
1 \to \GG_m \to X \to A \to 0,
$$
of a principally polarized abelian variety of dimension $g-1$.
Under the map $q$ this cycle is mapped to ${\Acal}_{g-1}$ in the
Satake compactification.
We denote by $\delta_g$ the cycle class in the sense
of the $\QQ$-classes, $[\Delta_g]_{\QQ}$,
of this codimension $g$ cycle
in the Chow group with rational coefficients 
of codimension $g$ cycles on  ${\Acal}_g^{\prime}$.
Note that for $g>1$ (resp.\ $g=1$)
the generic semi-abelian variety which is a trivial
extension by a rank $1$ torus has $4$ (resp.\ $2$)
automorphisms, so $\Delta_g$
is `counted with multiplicity $1/4$ (resp.\ $1/2$)'. 
We refer to \cite{Kr}, \cite{L-MB}
for cycle theory on stacks, but see also \cite{M1}, \cite{M3} and \cite{E-vdG}.
We now can formulate our result.
 
\begin{theorem}\label{main} In the Chow group $CH^g_{\QQ}({\Acal}_g^{\prime}\otimes k)$
of codimension $g$ cycles of
the moduli stack of rank $\leq 1$ degenerations
${\Acal}_g^{\prime}\otimes k$ we have the formula
$$\lambda_g= (-1)^g \zeta(1-2g) \, \delta_g,$$
where $\delta_g$ is the $\QQ$-class of the locus $\Delta_g$
of semi-stable abelian varieties which
are trivial extensions of an abelian variety of dimension $g-1$
with $\GG_m$.
\end{theorem}
\bigskip

Recall that $\zeta(1-2g)$ is a rational number and equals $-b_{2g}/2g$,
with $b_{2g}$ the $2g$th Bernoulli number.
\noindent
\begin{example}
We have $12\, \lambda_1 = \delta_1$,  $120\, \lambda_2=\delta_2$
and $252 \, \lambda_3=\delta_3$.
\end{example}

For $g=2$ and $g=3$ there is a canonical toroidal compactification
$\tilde{\Acal}_g$ of ${\Acal}_g$, the Delaunay-Voronoi compactification. In
\cite{vdG1} the second author obtained the following formulas for $\delta_g$: in
the rational Chow ring of $\tilde{\Acal}_g$ for $g=2$ and $g=3$
$$\delta_2= 120\lambda_2-\sigma_2, \qquad
\delta_3= 252\lambda_3-15 \lambda_1^2 \sigma_1+2\lambda_1\sigma_2,
$$
where $\sigma_i$ denotes a certain class of codimension $i$ lying in the
boundary.  Our formula gives the part that does not depend on the choice of
compactification.

Just as in the case of curves it is possible to introduce the tautological ring
for compactified moduli of abelian varities. One simply takes the subring of the
Chow ring generated by the $\lambda_i$ in toroidal compactification (for which
the Hodge bundle has been given a toroidal extension). This is easily seen to be
independent of the chosen compactification and its relations are easily
specified without making reference to a toroidal compactification. However, as
in the case of curves one would like to express natural loci as classes in the
tautological ring. This is somewhat problematic particularly when these classes
lie in the boundary as it is not even clear that this question is independent of
the toroidal compactification. We would like to suggest introducing instead the
\Definition{tautological module}, which by definition is the pushdown of the
tautological ring to the Satake compactification. Note that as the Satake
compactification is (highly) singular the tautological module is only a subspace
of the (rational) Chow group. We end this paper by giving some examples of how
to express the classes of natural loci as elements of the canonical module.
\end{section}
\begin{section}{The Proof of Theorem (1.1).}

In order to prove the theorem we may work on a level cover of the moduli space
${\Acal}_g^{\prime}$ for some level $n\geq 3$ prime to the characteristic of
the field $k$ and prove the corresponding
relation $\lambda_g= (-1)^g \zeta(1-2g) \, 
n \, \delta_g^{(n)}$ there. Here $\delta_g^{(n)}$ denotes the locus of
semi-abelian varieties with level-$n$ structure which are trivial extensions
of an $g-1$-dimensional abelian variety by a rank $1$ algebraic
torus. This has the
advantage that we can avoid the problems due to the existence of automorphisms.
In the proof we then have to employ an index $(n)$ for all objects.
Having said that we will carry out the computation by formally working in
level $1$ and assuming that the reader knows how to interpret our identities. 

In the computation we shall need a description of the space
${\Acal}_g^{\prime}$. We shall assume that the reader is familiar
with the construction of toroidal compactifications of ${\Acal}_g$.
It might help the reader to have a look at Mumford's paper \cite{M2},
where the moduli space of rank $1$ degenerations is used.
Using the natural map $q: {\Acal}_g^{\prime}\to
{\Acal}_g^*$ an \'etale cover of 
the divisor $B_g:= {\Acal}_g^{\prime} \backslash {\Acal}_g$
can be identified with the dual of the universal family 
$\hat{\Xcal}_{g-1} \to {\Acal}_{g-1} 
$ and using the principal polarization it can be identified with the universal
family ${\Xcal}_{g-1} \to {\Acal}_{g-1}$. The cycle $\Delta_g$
has as support 
the image of the zero section $s: {\Acal}_{g-1} \to \hat{\Xcal}_{g-1}$.

We recall how a point of $\hat{\Xcal}_{g-1}$ determines a semi-abelian variety.
If $Z$ is a principally polarized abelian variety of dimension $g-1$
with theta divisor $\Xi$ then the dual abelian variety $\hat{Z}$
classifies extensions
$$
1 \to \GG_m \to G \to Z \to 0
$$
of $Z$ by $\GG_m$.
Since the polarization defines a isomorphism $Z \to \hat{Z}$ we can associate
a semi-abelian variety to a point $z \in Z$. 
We may view this $\GG_m$-extension
as a $\GG_m$-bundle over $Z$ and we can take the corresponding $\PP^1$-bundle 
$\rho: \tilde{G} \to Z$. We now glue the $0$-section $\tilde{G}_0$ and 
the $\infty$-section $\tilde{G}_{\infty}$ over a translation by $z$ 
to get a non-normal variety $\bar{G}$.  Then $O(\tilde{G}_{\infty}+
\rho^{-1}(\Xi))$ descends to a line bundle $L$ on $\bar{G}$ with $h^0(L)=1$.
In this way we find a compactified semi-abelian variety canonically
associated to the pair $((Z,\Xi),z)$.

By doing this globally we see that
the moduli stack ${\Acal}_g^{\prime}$ comes with a universal semi-abelian
variety $\pi^{\prime}: {\Xcal}_g^{\prime}\to {\Acal}_g^{\prime}$ 
and a relative compactification $\bar{\pi}^{\prime}:
\bar{\Xcal}_g^{\prime} 
\to {\Acal}_g^{\prime}$ (i.e., this map $\bar{\pi}^{\prime}$ is proper).
To get it one
takes a smooth compactification $\pi: \tilde{\Xcal}_g \to \tilde{\Acal}_g$
as constructed in \cite{F-C} and restricts to $\pi^{-1}({\Acal}_g^{\prime})$.
The result then does not depend of the choice of  $\tilde{\Xcal}_g $.
The `universal' semi-abelian variety $G$ over the \'etale cover
$\hat{\Xcal} _{g-1}$ of $B_g \subset {\Acal}_g^{\prime}$ 
is the $\GG_m$-bundle obtained from 
the Poincar\'e bundle $P \to {\Xcal}_{g-1} \times {\hat {\Xcal} }_{g-1}$ 
by deleting the zero-section. 
We have the maps
$$
G=P -\{ (0)\} \to {\Xcal}_{g-1} \times_{{\Acal}_{g-1}} \hat{\Xcal}_{g-1}
{\buildrel q \over \to} {\Acal}_{g-1}.
$$

The fibre over $x\in B_g$ in the compactification $ \bar{\Xcal}_g^{\prime}$
of the universal family ${\Xcal}_g^{\prime}$ 
of semi-abelian varieties
is a compactification $\bar{G}$ of a $\GG_m$-bundle
$G$ over an abelian variety $X_{g-1}$ of dimension $g-1$ as constructed above. 
(In level $n\geq 3$ it is a compactification of a $\GG_m \times \ZZ/n\ZZ$-bundle
over an abelian variety $X_{g-1}$ 
of dimension $g-1$ and is a family of $n$-gons
over $X_{g-1}$.)
The points where $\bar{\pi}^{\prime}$
is not smooth are exactly the points of $\bar{G}-G$. 
So globally the locus where $\bar{\pi}^{\prime}$ is not smooth is the
codimension $2$ cycle $D$ in $\bar{\Xcal}_g^{\prime}$
obtained from gluing by a shift the $0$-section
and the $\infty$-section of the $\PP^1$-bundle associated to the Poincar\'e
bundle $P$ over ${\Xcal}_{g-1}\times_{{\Acal}_{g-1}} \hat{\Xcal}_{g-1}$.
We may identify an \'etale cover of the support of
$D$ with ${\Xcal}_{g-1} \times_{{\Acal}_{g-1}} {\hat {\Xcal}}_{g-1}$.

Our proof of Theorem (1.1) is based on an application of the 
Grothendieck-Riemann-Roch theorem to the structure sheaf 
on the universal semi-abelian variety over ${\Acal}_g^{\prime}$.
We start with a calculation on a smooth compactification  $\tilde{\Xcal}_g$
as constructed in \cite{F-C}  of  the universal
semi-abelian variety. We let 
$\pi\colon \tilde{\Xcal}_g \to \tilde{\Acal}_g$
be the natural morphism; if we restrict $\pi$ to ${\Acal}_g^{\prime}$ we get
$\bar{\pi}^{\prime}: \bar{\Xcal}_g^{\prime} \to {\Acal}_g^{\prime}$.

Applying the Grothendieck-Riemann-Roch theorem
to the  structure sheaf $O_{\tilde{\Xcal}_g}$  gives
in the Chow rings with rational coefficients 
$$
{\rm ch}(\pi_!O_{\tilde{\Xcal}_g})=\pi_*(e^{{\rm ch}(O_{\tilde{\Xcal}_g})} 
\, {\rm Td}^{\vee}(
\Omega^1_{\tilde{\Xcal}_g/\tilde{\Acal}_g}))=
\pi_*({\rm Td}^{\vee}(
\Omega^1_{\tilde{\Xcal}_g/\tilde{\Acal}_g})).$$
Here ${\rm Td}^{\vee}$ is the Todd class (that for a line bundle $L$ equals
$c_1(L)/(e^{c_1(L)}-1)$).
The relative cotangent sheaf fits in an exact sequence
$$
0 \to \Omega^1_{\tilde{\Xcal}_g/\tilde{\Acal}_g} \to \pi^*(\EE)
\to {\Fcal} \to 0
$$
with $\EE$ the Hodge bundle on $\tilde{\Acal}_g$
and ${\Fcal}$ a sheaf with support
where $\pi$ is not smooth.
Note that by \cite{F-C}, p.\ 195  we have
$$
\pi^*(\EE)= \Omega^1_{\tilde{\Xcal}_g}(\log)/
\pi^*(\Omega^1_{\tilde{\Acal}_g}(\log)),
$$
where $\log$ refers to logarithmic poles along the divisors at
infinity of $\tilde{\Xcal}_g$ and $\tilde{\Acal}_g$.

Substituting this in the Riemann-Roch formula we get
$$
{\rm ch}(\pi_!(O_{\tilde{\Xcal}_g}))=\pi_*(F) {\rm Td}^{\vee}(\EE)
$$
with $F := \mathrm{Td}^{\vee}({\Fcal})^{-1}$.
Since the cohomology of an abelian variety is the exterior algebra on $H^1$
the derived sheaf $\pi_!(O_{\tilde{\Xcal}_g})$ equals $\wedge^* \EE^{\vee}=
\sum_{i=0}^g (-1)^i \wedge^i \EE^{\vee}$. By the Borel-Serre formula \cite{B-S}, p.\
128 
we have ${\rm ch}(\wedge^* \EE^{\vee})= \lambda_g {\rm Td}(\EE)^{-1}$. Comparing
the terms of degree $\leq g$ in the resulting identity
$$
\lambda_g {\rm Td}(\EE)^{-1} = \pi_*(F){\rm Td}^{\vee}(\EE)
$$
yields  the result  of \cite{vdG1} :

\begin{proposition}
We have $\pi_*({\rm Td}^{\vee}({\Fcal})^{-1})=\pi_*(F)=\lambda_g$.
\end{proposition}

We now restrict to $\bar{\Xcal}_g^{\prime}$ and ${\Acal}_g^{\prime}$.
The sheaf ${\Fcal}$ has support on $D$. If $u$ is a fibre coordinate
on the $\GG_m$-bundle over the abelian scheme 
${\Xcal}_{g-1}\times_{{\Acal}_{g-1}}\hat{\Xcal}_{g-1}$ over $\hat{\Xcal}_{g-1}$
then a section of the pull back $\pi^*(\EE)$ of the Hodge bundle is given by
$du/u$. We now pull the section back 
to the $\PP^1$-bundle and take the residue along
the $0$-section and the $\infty$-section. This gives an isomorphism
of sheaves on ${\Acal}_g^{\prime}$ 
$${\Fcal} \cong O_{\tilde{D}},$$
where $\tilde{D}$ is the double \'etale cover of $D$
corresponding to choosing the branches $0$ and $\infty$ in the $\PP^1$-bundle.

The normal bundle to an \'etale cover of $D$ given by 
${\Xcal}_{g-1} \times_{{\Acal}_{g-1}} \hat{\Xcal}_{g-1}$ is
then $N= P \oplus \tau^*(P^{-1})$ with $P$ the Poincar\'e bundle and
$\tau$ the map from ${\Xcal}_{g-1}\times_{{\Acal}_{g-1}} \hat{\Xcal}_{g-1}$
to itself defining the translation by which we glue the $0$-section
and the $\infty$-section of the $\PP^1$-bundle corresponding to the Poincar\'e
bundle. On points $\tau$ is given by
$\tau(x, \hat x)= (x+\hat x, \hat x)$.
(We identify $\hat{\Xcal}$ with ${\Xcal}$ if needed.)  
We write $\alpha_1= c_1(P)$ and
$\alpha_2= c_1(\tau^*(P^{-1}))$ for the first Chern classes.  On the space of
rank $\leq 1$ degenerations $\bar{\Xcal}_g^{\prime}
$ we then can write $[D]= \alpha_1\alpha_2$.

Let $i\colon D \to \bar{\Xcal}_g^{\prime}$ be the inclusion.
Then if we write
$$
{\rm Td}^{\vee}(L)=\frac{c_1(L)}{ (e^{c_1(L)}-1)}=
\sum_{k=0}^{\infty} \frac{b_k}{ k!} (c_1(L))^k
$$
with $b_k$ the $k$-th Bernoulli number,
we have (cf.\  Mumford [10], p.\ 303):
$$
\bar{\pi}^{\prime}_*({\rm Td}^{\vee}(O_D^{-1} -1)) =
\bar{\pi}^{\prime}_*\left( 
\sum_{k=1}^{\infty} \frac{(-1)^kb_{2k}}{(2k)!}
i_*(\frac{ \alpha^{2k-1} + \alpha_2^{2k-1}}{ \alpha_1 + \alpha_2})\right).
\eqno(1)
$$
Consider now the Poincar\'e bundle $P$ on $X \times \hat{X}$ for an abelian
variety $X$ of dimension $g-1$ and dual abelian variety $\hat{X}$.
If $T$ is a line bundle on $X \cong \hat{X}$ 
that represents (locally in the \'etale topology)
the principal polarization of $\hat{X}_{g-1}$
then $P=m^*T \otimes p^*T^{-1} \otimes \hat{p}^*T^{-1}$.
Employing the notation $p$ and $\hat{p}$ for the projections
on $X$ and $\hat{X}$ we find

\begin{align} \notag
\tau^*(P^{-1})&= \tau^*((m^*T)^{-1}) 
\otimes \tau^*p^* T \otimes \tau^*\hat{p}^*T \\ \notag
&=\tau^*((m^*T)^{-1})\otimes m^*T \otimes \hat{p}^*T\\ \notag
\end{align}
since $p\tau=m$ and $\hat{p}\tau=\hat{p}$.
We get
$$
P \otimes \tau^*(P^{-1})\cong \tau^*(m^*T)^{-1}) \otimes (m^*T)^{\otimes 2}
\otimes p^*T^{-1}.
$$
Restriction to a fibre $X\times \hat{x}$ gives
$$
t_{-2\hat{x}}(T^{-1}) \otimes (t_x^*T)^{\otimes 2} \otimes T^{-1}
$$
and by the Theorem of the Square this is trivial on such a fibre 
$X \times \hat{x}$. This implies that on 
${\Xcal}_{g-1} \times_{{\Acal}_{g-1}}
\hat{\Xcal}_{g-1}$ we have
$$
c_1(N)= c_1(P\otimes \tau^*(P^{-1}))= \alpha_1+\alpha_2=\hat{p}^*(\beta)
$$
with $\beta$ a codimension $1$ class on ${\hat{\Xcal}}_{g-1}$.
In order to
determine $\beta$ we may restrict to the other fibre $0\times \hat{X}$.
Then $P_{| 0\times\hat{X}}$ is trivial and $\tau^*(P^{-1}){| 0\times \hat{X}}
$ is the pull back of $P^{-1}$ from the diagonal. But assuming as we
may that $T$ is symmetric we find that $P$ restricted
to the diagonal is $T^{\otimes 2}$. So as a result we find on $X
\times \hat{X}$ 
that $ \beta= -2c_1(T)$ on $\hat X$
and we get an identity on $X \times \hat{X}$
$$
N= P \oplus P^{-1}\otimes {\hat p}^*(T^{-2}).\eqno(2)
$$
We can consider this as a global identity on $D$ by considering
this as a definition of the line bundle $\hat{p}^*(T)$ 
on ${\Xcal}_{g-1}\times_{{\Acal}_{g-1}}
\hat{\Xcal}_{g-1}$. The
line bundle  $T$ restricts in each fibre $\hat{X}$ of ${p}$ to $O(\Theta)$
with $\Theta$ the theta divisor.
Developing the terms in (1) we get expressions of the form

\begin{align} \notag
\bar{\pi}^{\prime}_*(i_*(\alpha_1+\alpha_2)^r(\alpha_1\alpha_2)^s)&=
\bar{\pi}^{\prime}_*(i_*({\hat p}^*(\beta^r))(\alpha_1\alpha_2)^s)\notag
\\ &= j_*(\beta^r{\phi}_*(D^s)),\notag \\ \notag
\end{align}
where $\phi$ is the restriction to the boundary 
$\bar{\Xcal}^{\prime}_g -{\Xcal}_g^{\prime}$
of $\bar{\pi}^{\prime}$ and 
$j\colon B_g \to
{\Acal}_g^{\prime}$ is the inclusion of the boundary of ${\Acal}_g^{\prime}$.
Moreover, we use $\bar{\pi}^{\prime} i = \hat{p}$ 
and abuse the notation $D$ also
for the $\QQ$-class of $D$.

We claim that for dimension reasons the only surviving terms are of the form
$j_*(\beta^r){\phi}_*(D^{g-1})$. Indeed, the fibres of $\phi$
have dimension $g-1$.
Thus by Proposition 2.1 the only term in (1) that
can contribute to $\lambda_g$ is:
$$
\frac{(-1)^g b_{2g}}{ (2g)! }
\bar{\pi}^{\prime}_*i_*((-1)^{g-1}(2g-1)(\alpha_1\alpha_2)^{g-1}).
$$
So we need to compute $\bar{\pi}^{\prime}_*i_*(D^{g-1})= \pi_*(D^g)$.
The identity $\pi_*({\rm Td}^{\vee}({\Fcal})^{-1})=\lambda_g$ 
of Prop. 2.1 implies that
$$
\bar{\pi}^{\prime}_*(F)=\frac{(-1) b_{2g} }{ (2g)! }
\bar{\pi}^{\prime}_*i_*((2g-1)(\alpha_1\alpha_2)^{g-1})=\lambda_g.
$$
Represent the line bundle 
$P$ by the divisor $\Pi$ and the line bundle $T$ by a divisor $T$ (by abuse of
notation) on
${\Xcal}_{g-1} \times_{{\Acal}_{g-1}} \hat{\Xcal}_{g-1} $. Then we have by (2)
that $\alpha_1=\Pi$ and $\alpha_2=-\Pi-2T$, so
$$
(-1)^{g-1}\alpha_1^{g-1}\alpha_2^{g-1}
=\Pi^{2g-2} + \sum_{r=1}^{g-1} \left(\frac{g-1}{r}\right)
\Pi^{2g-2-r}{\hat p}^*(2T)^{r}.
$$
Now  apply  GRR to the  bundle $P\otimes {\hat p}^*(O(nT))$ on
${\Xcal}_{g-1} \times_{{\Acal}_{g-1}} \hat{\Xcal}_{g-1}$
and the morphism $\hat{p}$; it says
$$
{\rm ch}({\hat p}_!(P\otimes {\hat p}^*O(nT)))= {\hat p}_*(e^{\Pi}) \,
 \cdot e^{nT} \cdot {\rm Td}^{\vee}(q^*\EE_{g-1}), \eqno(3)
$$
where $q^*(\EE_{g-1})$ is the pull back to $\hat{\Xcal}_{g-1}$ of the Hodge
bundle $\EE_{g-1}$ on ${\Acal}_{g-1}$.
But ${\hat p}_!(P\otimes \hat{p}^*O(nT))$ 
is a sheaf with support (in  codimension
$g-1$) over the zero section $S$. Once again
by applying GRR, this time to the inclusion $S \to \hat{\Xcal}_{g-1}$,
(cf., \cite{M4}, p.\ 65) we see that by viewing
$${\hat p}_!(P\otimes \hat{p}^*O(nT))$$ as a derived sheaf on ${\Xcal}_{g-1}$
we get $c_i({\hat p}_!(P\otimes \hat{p}^*O(nT))=0$ for $i<g-1$ and
$$c_{g-1}({\hat p}_!(P\otimes
\hat{p}^*O(nT)))= (-1)^{g-2}(g-2)!\,  [S].
$$
It then follows by comparing codimensions $g$ classes which are
coefficients of the same powers of $n$ on both sides of (3)
that

$${\hat p}_*(\Pi^{2g-2-r})T^r  = 0\quad {\rm   if} \quad  r\neq 0$$
and
$${\hat p}_*(\Pi^{2g-2-r})T^r=(-1)^{g-1}(2g-2)!\, [S] \quad
{\rm if } \quad   r=0.$$
So we find
$$
\pi_*(D^g)=j_*({\hat p}_*(\Pi^{2g-2}))= (-1)^{g-1} (2g-2)![\Delta_g],
$$
where $\Delta_g$ is the zero section of ${\Xcal}_{g-1} \to {\Acal}_{g-1}$.
Interpreting this identity in the right way means taking the $\QQ$-class
of $\Delta_g$.
We thus get

\begin{align} \notag
\lambda_g & = \frac{(-1) b_{2g} }{(2g)! }
\bar{\pi}^{\prime}_*i_*((2g-1)(\alpha_1\alpha_2)^{g-1}) \\ \notag
& = \frac{(-1) b_{2g} } { (2g)! }(2g-1) (2g-2)! (-1)^{g-1}[\Delta_g]_{\QQ} \\ \notag
&= (-1)^g \zeta(1-2g)[\Delta_g]_{\QQ}\\ \notag
\end{align}
as required. This concludes the proof of the theorem.
\end{section}
\begin{section}{The tautological module}

Recall that the Hodge bundle extends to a toroidal compactification
$\tilde{A}_g$ as the dual of the Lie algebra of the semi-abelian variety that is
supposed to exist over $\tilde{A}_g$. Recall also that the subring of
$CH^*_\QQ(\tilde{A}_g)$ generated by the Chern classes 
$\lambda_i$ of this extension
is independent of the choice of toroidal compactification. Indeed, 
the relation 
$$
(1+\lambda_1+\ldots + \lambda_g)
(1-\lambda_1+\lambda_2- \ldots +(-1)^g\lambda_g)=1
$$
(see \cite{vdG1}, \cite{E-V}) that always exists between the $\lambda_i$
suffices to show that this ring is a Gorenstein algebra with socle in the top
degree, $g(g+1)/2$, and is hence determined by the evaluation map in that degree
which is independent of the choice of $\tilde{A}_g$. It is still somewhat
unpleasant that this \emph{tautological subring} is a ring of classes living on
different toroidal compactifications. We would like to suggest that the pushdown
of classes in this this ring to the Satake compactification should also be of
interest. We begin by showing that it is independent of the choice of toroidal
compactification.
\begin{proposition}
Let $\tilde{A}_g$ be a toroidal compactification with
\map{q}{\tilde{A}_g}{\Acal^*_g} the canonical map to the Satake compactification
and let $\alpha$ be a subset of $\{1,2,\dots,g\}$. Then the class
$\ell_\alpha := q_*(\lambda_\alpha) \in CH^*_\QQ(\Acal^*_g)$ is independent of
$\tilde{A}_g$, where $\lambda_\alpha = \prod_{i \in \alpha}\lambda_i$.
\end{proposition}
\begin{proof}
Any two toroidal compactifications have a common refinement \cite[p.~97--98, i),
iii)]{F-C} so we may assume that one compactification is a refinement of the
other and then it is clear,
as the the $\lambda_i$ are compatible with pullback,
pulling and pushing down is the identity for a birational map,
and pushing down is transitive.
\end{proof}
Of particular interest is of course to which extent the cycles of natural
subvarieties of $\Acal^*_g$ can be expressed as linear combinations of the
$q_*(\lambda_\alpha)$. The class $\ell_{\{1\}}$ is actually the class of 
the natural ample line bundle on $\Acal^*_g$ 
so that the linear subspace sections represent
the classes $q_*(\lambda_1^i)$ which are (explicit) 
such linear combinations and
generally if the class of subvariety has such a representation then any
hyperplane section of it does too.

In positive characteristic a less trivial example can be found in
\cite{vdG1}. Consider the closed algebraic subset $V_0$ of ${\Acal}_g\otimes
\FF_p$ of all abelian varieties with $p$-rank zero. By Koblitz (see [K]) we know
that this is a pure codimension $g$ cycle on ${\Acal}_g\otimes \FF_p$. It is a
complete cycle since abelian varieties of $p$-rank $0$ cannot degenerate. Any
complete subvariety of ${\Acal}_g$ has codimension at least $g$ in ${\Acal}_g$,
see \cite{vdG1} and \cite{O}.
\begin{theorem} 
The $\QQ$-cycle class of $V_0$ in $CH^g_{\QQ}(\Acal^*_g\otimes \FF_p)$ is
given by the formula $[V_0]_{\QQ}=(p-1)(p^2-1)\cdots (p^g-1)\, \ell_{\{g\}}$.
\end{theorem}
\begin{remark}
i) Recently, Keel and Sadun proved that there is no complete subvariety of
codimension $g$ in ${\Acal}_g \otimes {\CC}$ for $g\geq 3$, cf.~\cite{K-S}.

ii) In positive characteristic there are many other natural subvarieties whose
classes will lie in the tautological module.
\end{remark}
We shall now show that the our main result can be used to express the top and
next to the top boundary component as tautological classes.
\begin{theorem} 
i) In the group $CH_{\QQ}^g({\Acal}_g^*)$
we have $[\Acal^*_{g-1}]=(-1)^g/\zeta(1-2g)\ell_{\{g\}}$.

ii) In the group $CH_{\QQ}^{2g-1}({\Acal}_g^*)$ 
we have $[\Acal^*_{g-2}]=(1/\zeta(1-2g)\zeta(3-2g))
 \ell_{\{g-1,g\}}$.
\end{theorem}
\begin{proof}
For the first part we note that by the excision exact sequence
$CH_{\QQ}^g(\Acal^*_{g-2}) \to CH_{\QQ}^g({\Acal}_g^*) \to CH_{\QQ}^g
({\Acal}_g\setminus\Acal^*_{g-2}) \to 0$ we may prove it in $\Acal^*_g\setminus
\Acal^*_{g-2}$. We have a proper map $\Acal'_g \to \Acal^*_g\setminus
\Acal^*_{g-2}$ and hence the formula follows from theorem \ref{main} by pushing
down the main formula to $\Acal^*_g$. As for the second we have that on
$\Acal'_g$ $\lambda_g\lambda_{g-1}= (-1)^g/\zeta(1-2g)\lambda_{g-1}\delta_g$ and
the support of $\delta_g$ maps finitely to $\Acal_{g-1}$ and the restriction of
$\lambda_{g-1}$ to it corresponds to $\lambda_{g-1}$ on $\Acal_{g-1}$ which is
zero. Hence, $\lambda_g\lambda_{g-1}$ is zero on $\Acal'_g$ and for dimension
reasons and excision is a multiple of $[\Acal^*_{g-2}]$. That multiple can be
determined by intersecting with $\lambda_1^{(g-1)(g-2)/2}$. Using the 
intersection numbers in \cite{vdG1}, p.\ 72 one sees that
$$
\lambda_g\lambda_{g-1} \lambda_1^{(g-1)(g-2)/2}[\tilde{\Acal}_g]
 = (1/\zeta(1-2g)\zeta(3-2g))
\lambda_1^{(g-1)(g-2)/2}[\tilde{\Acal}_{g-2}]
$$ 
and then one uses $\lambda_1^{(g-1)(g-2)/2}[\tilde{\Acal}_{g-2}]=
\lambda_1^{(g-1)(g-2)/2}[{\Acal}_{g-2}^*]$.
\end{proof}

The statements in the preceding theorem suggest an immediate
generalization. To lend some credibility to such a generalization
we sketch a proof of such a generalization in positive characteristic.

\begin{theorem} In characteristic $p>0$ we have in
$CH^{d}_{\QQ}({\Acal}_g^*)$ with $d=g(g+1)/2-(g-i)(g+1-i)/2$ the relation
$$
[A_{g-i}^*]= (-1)^i \frac{1}{\prod_{j=1}^i \zeta(2j-1-2g)} \ell_{\{g-i+1,
\ldots,g\}}.
$$
\end{theorem}
\begin{proof} The idea of the proof is to use the fact that
the cycle class of the locus $V^{(g)}_f$ in $\tilde{\Acal}_g$ of
semi-abelian varieties of $p$-rank $\leq f$ is a non-zero multiple
of $\lambda_{g-f}$, cf.\ \cite{vdG1}, and the fact that a semi-abelian
variety of $p$-rank $\leq f$ has torus rank $\leq f$.

We use a toroidal compactification $\tilde{\Acal}_g$ of Faltings-Chai type.
The closure $\bar{\Delta}_g$ in $\tilde{\Acal}_g$
of the locus $\Delta_g$ has a cover which is
a toroidal compactification ${\Acal}_{g-1}^{T}$ 
of ${\Acal}_{g-1}$, but not necessarily smooth.
But there is a smooth toroidal compactification $\tilde{\Acal}_{g-1}$
and a morphism $\tilde{\Acal}_{g-1} \to {\Acal}_{g-1}^{T}$ with the
property that the pull-back of the universal semi-abelian 
variety $\tilde{\Xcal}_g$ to $\tilde{\Acal}_{g-1}$ is a product of a
universal semi-abelian variety $\tilde{\Xcal}_{g-1}$ with a torus.
Then this can be repeated: the smooth toroidal compactification
$\tilde{\Acal}_{g-1}$ is a compactification of the canonical
partial compactification ${\Acal}_{g-1}^{\prime}$ and it contains a locus
$\Delta_{g-1}$ corresponding to trivial extensions of $(g-2)$-dimensional
abelian varieties with $\GG_m$. 
We can then consider the closure $\bar{\Delta}_{g-1}$ and there is a smooth
toroidal compactification $\tilde{\Acal}_{g-2}$ mapping to 
$\bar{\Delta}_{g-1}$ with a similar property.

>From our theorem \ref{main} know that we have a relation in 
$CH_{\QQ}^*(\tilde{\Acal}_g)$
$$
\lambda_g \doteq [\bar{\Delta}_g]+ \sigma_g,
$$
where $\sigma_g$ is a class with support on $q^{-1}({\Acal}_{g-2}^*)$
and where $\doteq$ means equality up to a
non-zero multiplicative factor which is a rational number.
The relation $[V_1^{(g)}] \doteq \lambda_{g-1}$ in $CH^{g-1}_{\QQ}(
\tilde{\Acal}_g)$ gives
$$
\lambda_g \lambda_{g-1} \doteq [V_1^{(g)}] \cdot [\bar{\Delta}_g] +
[V_1^{(g)}]\cdot \sigma
\doteq [V_1^{(g)}]\cdot [\bar{\Delta}_g],\\ 
$$
since $V_1^{(g)}$ does not intersect the torus rank $\geq 2$ locus $q^{-1}(
{\Acal}^*_{g-2})$.
Now the pull back of the intersection of $V_1^{(g)}$ with $\bar{\Delta}_g$
to $\tilde{\Acal}_{g-1}$ is exactly the $p$-rank zero locus $V_0^{(g-1)}$
on ${\Acal}_{g-1}$.
But we know that 
$$
[V_0^{(g-1)}]\doteq \lambda_{g-1} \doteq \delta_{g-1}
$$
on ${\Acal}_{g-1}^{\prime}$, hence on the cover $\tilde{\Acal}_{g-1}$
of $\bar{\Delta}_{g}$ we get the relation
$$
[V_1^{(g-1)}] \cdot \bar{\Delta}_{g} \doteq \delta_{g-1} + \sigma_{g-1}
$$
with $\sigma_{g-1}$ a class with support on $q^{-1}({\Acal}_{g-3}^*)$.
Now we use the relation
$
\lambda_{g-2} \doteq [V_2^{(g)}]
$
and get
\begin{align} \notag
\lambda_g \lambda_{g-1} \lambda_{g-2} & \doteq 
[V_2^{(g)}] \cdot \bar{\Delta}_{g-1} + [V_2^{(g)}] \cdot \sigma_{g-1} \\ \notag
&\doteq [V_2^{(g)}] \cdot \bar{\Delta}_{g-1}\\ \notag
\end{align}
and the pull back of the intersection of $V_2^{(g)}$ with $\bar{\Delta}_{g-1}$
to $\tilde{\Acal}_{g-2}$ is exactly $V_0^{(g-2)}$.
But using again a relation 
$[V_0^{(g-2)}]\doteq \lambda_{g-2} \doteq \delta_{g-1}$ 
on ${\Acal}_{g-2}^{\prime}$ we see that a non-zero multiple
$\lambda_g \lambda_{g-1} \lambda_{g-2}$ is represented by the cycle 
$\delta_{g-1}$ 
on $\tilde{\Acal}_g-q^{-1}({\Acal}_{g-3}^*)$. Arguing as in the proof of the
preceding theorem we see that 
$q_*(\lambda_g \lambda_{g-1} \lambda_{g-2})$ 
is a non-zero multiple of $[A_{g-2}^*]$
To determine the multiple we intersect again with the appropriate power
of $\lambda_1$, cf.\ \cite{vdG1}. Proceeding in this way by induction one deduces the theorem.
\end{proof}

By a simple argument we can also show that another class is also in the
tautological module.
\begin{proposition}
The cycle class $[B_g^*]$ of the boundary is the same in the Chow group
$CH_{\QQ}^g({\Acal}_g^*)$ as a multiple of the $\QQ$-class of the locus of
products $X \times E$ of principally polarized abelian varieties of dimension
$g-1$ with a fixed elliptic curve~$E$.
\end{proposition}
\begin{proof}
For $g=1,2$ see \cite{vdG2}. Consider (for $g>2$) the space ${\Acal}_{g-1,1}$ of
products of a principally polarized abelian variety of dimension $g-1$ and an
elliptic curve. It is the image of ${\Acal}_{g-1} \times {\Acal}_1$ in
${\Acal}_g$ under a morphism to ${\Acal}_g$ which can be extended to a morphism
${\Acal}_{g-1}^* \times {\Acal}_1^* \to {\Acal}_g^*$.  Since an \'etale cover of
${\Acal}_1$ is the affine $j$-line we find a rational equivalence between the
cycle class of a fibre ${\Acal}^*_{g-1} \times \{ j \}$ with $j$ a fixed point
on the $j$-line and a multiple of the fundamental class of the boundary $B_g^*$.
\end{proof}

\end{section}

\newcommand\eprint[1]{Eprint:~\texttt{#1}}

%
%


\begin{thebibliography}{9999}

\bibitem{B-S} A.\ Borel, J-P. Serre: Le th\'eor\`eme de Riemann-Roch (d'apr\`es
Grothendieck). {\sl Bull.\ Soc.\ Math.\ France \bf 86} (1958), 97--136.

\bibitem{E-vdG} T.\ Ekedahl, G.\ van der Geer: The order of the top
Chern class of the Hodge bundle on the moduli space of abelian varieties.
To appear in Acta Mathematica. (math.AG/0302291)

\bibitem{E-V} H.\ Esnault, E.\ Viehweg: Chern classes of Gauss-Manin 
bundles of weight 1 vanish.  {\sl $K$-Theory  \bf 26}  (2002), 287--305.

\bibitem{F-C} G.\ Faltings, C-L.\ Chai: Degeneration of abelian varieties.
Ergebnisse der Math. 22. Springer Verlag 1990.

\bibitem{Fu} W.\ Fulton: Intersection theory. Ergebnisse der Math. 2. Springer
Verlag 1984.

\bibitem{vdG1} G.\ van der Geer: Cycles on the moduli
space of abelian varieties. In {\sl Moduli of Curves and Abelian
Varieties}, The Dutch Intercity Seminar on Moduli.
Eds. C.\ Faber and E.\ Looijenga. Vieweg Verlag 1999. Braunschweig.

\bibitem{vdG2} G.\ van der Geer: The Chow ring of the moduli space of abelian threefolds.
{\sl J.\ of Alg.\ Geometry \bf 7} (1998), 753--770.

\bibitem{K-S} S.\ Keel, L.\ Sadun: Oort's conjecture for $A_g \otimes {\CC}$.
\eprint{math.AG/0204229}

\bibitem{K} N.\ Koblitz: $p$-adic variation of the zeta function over families
of varieties defined over finite fields. {\sl Compositio Math.\ \bf 31}, (1975),
119--218.

\bibitem{Kr} A.~Kresch, \emph{Cycle groups for {A}rtin stacks}, Invent. Math. \textbf{138}
  (1999), no.~3, 495--536, \eprint{math.AG/9810166}.

\bibitem{L-MB} G.\ Laumon, L.\ Moret-Bailly, Champs alg\'ebriques,
Springer-Verlag, 2000.

\bibitem{M1} D.\ Mumford: Picard groups of moduli problems. In {\sl
Arithmetical Algebraic Geometry}. Ed.\ O.F.G.\ Schilling.
Harper and Row, 1965, p.\ 33-81. 

\bibitem{M2} D.\ Mumford: On the Kodaira dimension of the Siegel modular
variety. I2: Algebraic Geometry-Open Problems. SLNM 997, 348-375.

\bibitem{M3} D.\ Mumford: Towards an enumerative geometry of the
moduli space of curves. In: {\sl 
 Arithmetic and Geometry}, Vol.\ II. Birkh\"auser Verlag, Boston, MA,
1983, p.\ 271--328. 

\bibitem{M4} D.\ Mumford: 
Stability of projective varieties. Monographie de L'Enseignement Math.\ No.\ 24.Geneva 1977.

\bibitem{O} F.\ Oort: Complete subvarieties of moduli spaces. In:
{\sl Abelian
Varieties} (W.\ Barth, K.\ Hulek, H.\ Lange, eds.), de Gruyter Verlag,
Berlin, 1995, p. 225--235.


\end{thebibliography}
\end{document}